\documentclass[11pt, reqno]{amsart}

\usepackage{amssymb,geometry,stmaryrd,color}
\usepackage[T1]{fontenc}

\evensidemargin\oddsidemargin

\begin{document}

\baselineskip=18pt
\setcounter{page}{1}
    
\newtheorem{Prop}{Proposition}
\newtheorem{Lemm}{Lemma}
\newtheorem{Rem}{Remark}
\newtheorem{Coro}{Corollary}
\newtheorem{Propo}{Theorem}

\def\a{\alpha}
\def\b{\beta}
\def\B{{\bf B}} 
\def\CC{{\mathbb{C}}} 
\def\cN{{\mathcal{N}}} 
\def\cB{{\mathcal{B}}} 
\def\cI{{\mathcal{I}}} 
\def\cS{{\mathcal{S}}}
\def\UU{{\mathcal{U}}}
\def\ca{c_{\a}}
\def\ka{\kappa_{\a}}
\def\coa{c_{\a, 0}}
\def\cua{c_{\a, u}}
\def\cL{{\mathcal{L}}} 
\def\cW{{\mathcal{W}}} 
\def\Ea{E_\a}
\def\Eab{E_{\a,\b}}
\def\eps{{\varepsilon}} 
\def\esp{{\mathbb{E}}} 
\def\g{{\gamma}}
\def\Ga{{\Gamma}} 
\def\G{{\bf \Gamma}} 
\def\hag{{h_{\a,\g}}}
\def\thag{{{\hat h}_{\a,\g}}}
\def\K{{\bf K}}
\def\HH{{\bf H}}
\def\ii{{\rm i}}
\def\e{{\rm e}}
\def\L{{\bf L}}
\def\lbd{\lambda}
\def\sg{\sigma}
\def\lacc{\left\{}
\def\lcr{\left[}
\def\lpa{\left(}
\def\lva{\left|}
\def\M{{\bf M}}
\def\Ma{\M_\a}
\def\Mab{\M_{\a,\b}}
\def\Mabe{\M_{\a,\b, \varepsilon}}
\def\NN{{\mathbb{N}}} 
\def\pb{{\mathbb{P}}} 
\def\tpab{\hat{\varphi}_{a,b}} 
\def\tpa{\hat{\psi}_{\a}}
\def\tppa{\tilde{\psi}_{\a}} 
\def\tva{\hat{\varphi}_{\a}} 
\def\rl{{\mathbb{R}}}
\def\racc{\right\}}
\def\rpa{\right)}
\def\rcr{\right]}
\def\rva{\right|}
\def\prost{{\succ_{\! st}}}
\def\I{{\bf I}}
\def\X{{\bf X}}
\def\tX{{\hat \X}}
\def\Xag{{\X_{\a,\g}}}
\def\tXag{{{\hat \X}_{\a,\g}}}
\def\Y{{\bf Y}}
\def\tY{{\hat \Y}}
\def\Yag{{\Y_{\a,\g}}}
\def\tYag{{{\hat \Y}_{\a,\g}}}
\def\Z{{\bf Z}}
\def\Za{\Z_{\a}}
\def\Xab{\X_{\a,\b}}
\def\Yab{\Y_{\a,\b}}
\def\Yabe{\Y_{\a,\b, \varepsilon}}
\def\XX{{\mathcal X}}
\def\Y{{\bf Y}}
\def\J{{\bf J}}
\def\V{{\bf V}_\a}
\def\Un{{\bf 1}}
\def\ZZ{{\mathbb{Z}}}
\def\T{{\bf T}}
\def\AA{{\mathcal A}}
\def\hAA{{\hat \AA}}
\def\hL{{\hat L}}
\def\hT{{\hat T}}

\def\claw{\stackrel{d}{\longrightarrow}}
\def\elaw{\stackrel{d}{=}}
\def\qed{\hfill$\square$}

\newcommand*\pFqskip{8mu}
\catcode`,\active
\newcommand*\pFq{\begingroup
        \catcode`\,\active
        \def ,{\mskip\pFqskip\relax}%
        \dopFq
}
\catcode`\,12
\def\dopFq#1#2#3#4#5{%
        {}_{#1}F_{#2}\biggl[\genfrac..{0pt}{}{#3}{#4};#5\biggr]%
        \endgroup
}

\title{A note on the $\alpha-$Sun distribution}

\author[Thomas Simon]{Thomas Simon}

\address{Laboratoire Paul Painlev\'e, Universit\'e de Lille, 59000 Lille. {\em Email}: {\tt thomas.simon@univ-lille.fr}}



\begin{abstract} We investigate the analytical properties of the $\a-$Sun random variable, which arises from the domain of attraction of certain storage models involving a maximum and a sum. In the Fréchet case we show that this random variable is infinitely divisible,  and we give the exact behaviour of the density at zero. In the Weibull case we give the exact behaviour of the density at infinity, and we show that the behaviour at zero is neither polynomial nor exponential. This answers the open questions in the recent paper \cite{GW}.
\end{abstract}

\maketitle

\section{Introduction}
 
The $\a-$Sun random variable is a positive random variable $\Xag$  parametrized by $\a\in (0,1)$ and $\g > 0,$ whose smooth density $h_{\a, \g}$ solves the integro-differential equation
\begin{equation}
\label{ID}
\hag(x)\; =\; \frac{\g}{x}\int_0^x \frac{\hag(u)}{(x - \a u)^\g}\, du
\end{equation}
on $(0,\infty).$ This random variable was introduced in \cite{GH} as the renormalized limit of the sequence $\{Y_n,\, n\ge 0\}$ defined by
$Y_0 = X_0$ and
\begin{equation}
\label{Mod}
Y_n = \max \{Y_{n-1}, \alpha Y_{n-1} +X_n\}
\end{equation}
for $n\ge 1,$ where $\{X_n, \, n\ge 0\}$ is a given i.i.d. sequence belonging to the domain of attraction of a Fréchet distribution with index $\g.$ The boundary cases $\a = 0$ and $\a = 1,$ which will not be considered here, give a Fréchet respectively a positive stable random variable, as can be proved either from \eqref{ID} or from the recurrence \eqref{Mod} and the classical limit theorems. The $\a-$Sun random variable is hence a kind of interpolation between the Fréchet and the positive stable random variable. We refer to the recent paper \cite{GW} and the references therein for more detail on the $\a-$Sun random variable and the related storage models. 

The paper \cite{GW} undertakes an analytical study of the density $\hag$ starting from the following equation on the Mellin transform
\begin{equation}
\label{Mell}
\esp[\X_{\a,\g}^{s-1}]\; = \; \frac{\g\, \esp[\X_{\a,\g}^{s-\g-1}]}{(1+\g -s)}\; \pFq{2}{1}{\g,,1+\g-s}{2+\g-s}{\a}
\end{equation}
for all $s < \g +1,$ which is obtained from a direct integration of \eqref{ID} and Fubini's theorem - see (3.2) in \cite{GW}. This equation can be solved, and gives a Mellin-Barnes integral representation for $\hag$ which is not classical - see (3.38) and (3.39) in \cite{GW}. In particular the density function $\hag$, which is real analytic on $(0,\infty),$ cannot seem to be included in the large class of Fox $H-$functions. The authors also argue that $\Xag$ should be infinitely divisible (ID) as a renormalized limit, and ask whether this can be seen transparently from \eqref{Mell} - see the discussion in \cite{GW}. One purpose of this note is to show that this is indeed the case, and that the law of $\Xag$ actually belongs to a certain subclass of ID distributions, which we now define. The law of a positive random variable $X$ is called a generalized Gamma convolution (GGC) if $X$ admits an integral representation
$$X\; \elaw\; \int_0^\infty a(t)\, d\G_t$$
where $\{a(t), \; t\ge 0\}$ is a suitably integrable deterministic function and $\{\G_t, \; t\ge 0\}$ stands for the usual Gamma subordinator, here and throughout. Equivalently, the random variable is ID with a Lévy measure having a density $f$ such that the function $xf(x)$ is completely monotone on $(0,\infty)$. We refer e.g. to chapters 3 and 4 in \cite{Bo} for more material on this family of ID distributions.
  
\begin{Propo} 
\label{ggc}
For every $\a\in (0,1)$ and $\g > 0,$ the law of $\X_{\a,\g}$ is a {\em GGC}.
\end{Propo}

The proof of this result relies on the representation of $\Yag = \X^{-\g}_{\a,\g}$ as the terminal value of some multiplicative martingale involving powers of cut-off beta random variables. The negative powers of order $-1/\g$ of the latter turn out to all have a GGC distribution, and we can then appeal to a powerful result by Bondesson \cite{Bond} stating that the GGC class is stable under independent multiplication. The same kind of argument was already used in \cite{JSW} for the generalized stable densities, which satisfy an integro-differential equation close to \eqref{ID}. 

\medskip   

The exact behaviour of $h_{\a,\g}$ at infinity is easily derived from \eqref{ID} and reads 
$$\hag(x)\; =\; \frac{\g}{x^{\g+1} }\int_0^x \frac{\hag(u)}{(1 - \a ux^{-1})^\g}\, du\; \sim\; \g\, x^{-\g-1},$$
by monotone convergence and the fact that $\hag$ is a density function. This behaviour at infinity is that of the Fréchet random variable $\L^{-\frac{1}{\g}},$ where $\L$ stands for the unit exponential random variable. The possibility of a complete asymptotic expansion at infinity with explicit coefficients is discussed in \cite{GW} via the Mellin-Barnes representation - see (4.1) therein. The more subtle behaviour of $\hag$ at zero is however left there as an open problem - see the end of Section 4 in \cite{GW}. 

\begin{Propo} 
\label{asym}
For every $\a\in (0,1)$ and $\g > 0,$ one has
$$\hag(x)\; \sim\; c_{\a,\g}\, x^{-\frac{\g}{1-\a} -1}\, e^{-((1-\a) x)^{-\g}}$$
as $x\to 0,$ where 
$$ c_{\a,\g}\; =\; \g (1-\a)^{\frac{\g}{\a-1}} e^{\frac{\a\psi(1)}{\a-1}} \prod_{k=1}^\infty \, \frac{e^{\frac{\a}{(\a-1)k}}}{G_{\a,\g}(k)}\,\in\, (0,\infty)\qquad\mbox{with}\qquad G_{\a,\g}(k)\, =\, \pFq{2}{1}{\g,,1}{1+k\g}{\frac{\a}{\a-1}}.$$
\end{Propo}
In the above expression of $c_{\a,\g},$ the parameter $\psi(1)$ is the negative of Euler's constant, written this way in order to avoid confusion with the notation borrowed from \cite{GW} for the parameter $\g$, and the positive finite character of the infinite product defining $c_{\a,\g}$ is a consequence of the easily established asymptotic
$$G_{\a,\g}(k)\; =\; 1 \, +\, \frac{\a}{k(\a-1)}\, +\, O\lpa k^{-2}\rpa,\qquad k\to\infty.$$  
Observe that Theorem \ref{asym} shows that the behaviour of $\hag$ at zero is that of a generalized Fréchet random variable $(1-\a)^{-1} \G_{(1-\a)^{-1}}^{-\frac{1}{\g}},$ and that when $\a \to 0$ we recognize the density of the Fréchet random variable $\L^{-\frac{1}{\g}}.$ The proof of Theorem 2 is also inspired from \cite{JSW}, where the asymptotics at zero of the so-called generalized stable densities were studied via the underlying perpetuity of a  spectrally negative Lévy process - see Section 2.3 therein. The same argument could actually be applied here - see Remark \ref{SNLP} below, but we had rather exhibit the perpetuity of a subordinator associated to $\Yag,$ in order to use the recent asymptotic studies of \cite{H, MS} and show that they can be made completely explicit in the present context.

\medskip

The renormalized limit of \eqref{Mod} where the data sequence $\{X_n,\; n\ge 0\}$ belongs to the domain of attraction of a Weibull or a Gumbel distribution, has been studied in \cite{GHB}. The limit law of \eqref{Mod} in the Gumbel case has an explicit density on $\rl$ given by Corollary 2 in \cite{GHB}:
$$\frac{\exp\lpa -x - e^{-(1-\a)x}\rpa}{\Ga(1+(1-\a)^{-1})},$$
which is that of the random variable $-(1-\a)^{-1} \log\G_{(1-\a)^{-1}}.$ Observe in passing by Example 7.2.3 in \cite{Bond} that the latter random variable has a so-called extended GGC distribution. In the Weibull case, the limit law of \eqref{Mod} is that of a negative random variable $-\tXag,$ where $\tXag$ has a density $\thag$ solving the integro-differential equation 
\begin{equation}
\label{IDB}
\thag(x)\; =\; \frac{\g}{x}\int_x^{\frac{x}{\a}} (x - \a u)^\g\,\thag(u)\, du
\end{equation}
on $(0,\infty)$ - see Corollary 1 in \cite{GHB}. The latter equation is similar to \eqref{ID} but more difficult to analyze because both bounds in the integral depend on $x,$ except in the boundary case $\a = 0$ where the unique solution to \eqref{IDB} is the expected Weibull density ${\hat h}_{0,\g} (x) = \g\, x^{\g -1}\, e^{-x^\g}.$ The basic analytical properties of the density function $\thag$ do not seem to have been studied as yet - see the introduction of \cite{GW}. We show the following.

\begin{Propo}
\label{asymb}
For every $\a\in (0,1)$ and $\g > 0,$ the following holds.

\medskip

{\em (a)} The density $x\mapsto \thag(x)$ is strictly unimodal on $(0,\infty).$

\medskip

{\em (b)} For every $a > 0,$ one has $\lim_{x\to 0} x^{-a} \thag(x) \, =\, \lim_{x\to 0} x^a\, \log\thag(x) \, =\, 0.$ 

\medskip

{\em (c)} There exists $c\in (0,\infty)$ such that $\thag(x)\,\sim\, c\, x^{\frac{\g}{1-\a} -1}\, e^{-((1-\a) x)^{\g}}$ as $x\to \infty.$ 

\end{Propo}

Observe in particular that contrary to $\Xag,$ the random variable $\tXag$ is not infinitely divisible for $\g > 1$ by the well-known criterion on thin tails at infinity. One might conjecture that the law of $\tXag$ is a GGC for all $\a\in (0,1)$ and $\g\le 1,$ as is the case in the Weibull case $\a = 0$ - see Example 4.3.4 in \cite{Bo}. Unfortunately the terminal value property of Theorem \ref{ggc}, which holds here as well and helps prove Part (c), seems less useful for ID properties because the involved random variables have a support which is bounded away from zero and infinity for all $\a\in (0,1)$ - see Remark \ref{Mouss} (c) below. The asymptotics of Theorem \ref{asymb} show that $\tXag$ has integer moments of arbitrary order, positive or negative, and it is a natural question whether these integer moments characterize the law of $\tXag,$ in other words whether $\tXag$ or $\tX_{\a,\g}^{-1}$ is $M-$det. This question was addressed in Proposition 4.1 of our previous paper \cite{MB}, in the framework of generalized stable random variables. 

\begin{Propo}
\label{Momus}
For every $\a\in (0,1), \g> 0$ and $t\neq 0,$ one has
$$\tX^{t}_{\a,\g} \;\,\mbox{is {\em $M$-det}}\quad\Longleftrightarrow\quad 0\,<\, t\,\le\, 2\g\quad\Longleftrightarrow\quad \X^{-t}_{\a,\g}\;\, \mbox{is {\em $M$-det}.}$$
\end{Propo}

The reason why we consider power transformations in the above result is partly methodological, and the first equivalence gives actually the main argument for Theorem \ref{asymb} (b). For $\tXag,$ the proof of the characterization relies on the crucial observation that the function $\thag(e^x)$ is log-concave on $\rl$, which by Lin's condition implies that Krein's integral criterion for moment-indeterminacy is also a necessary condition. The log-concavity also immediately implies Theorem \ref{asymb} (a). For $\Xag,$ the log-concavity argument is no more valid and we use the alternative property that $\log\Xag$ is self-decomposable, which implies an extended Lin's property as observed in \cite{MB}.  

\medskip

This note is organized as follows. In the next section we show the four above results, in their order except for Theorem \ref{Momus} which is proved before Theorem \ref{asymb} (a) and (b). In the last section we comment on an identity between the perpetuity of a general  subordinator and the terminal value of a multiplicative martingale, which is the main theme of the present paper and which we further illustrate with an example taken from the recent paper \cite{Ber}. 

\section{Proofs} 

\subsection{Proof of Theorem \ref{ggc}} We consider the negative power transformation $\Yag = \X_{\a,\g}^{-\g},$ whose positive integer moments are expressed as
\begin{equation}
\label{Mom}
\esp[\Y_{\a, \g}^n]\; =\; \prod_{k=1}^n \, \frac{k}{F_{\a,\g}(k)}\qquad\mbox{with}\qquad F_{\a,\g}(k)\, =\, \pFq{2}{1}{\g,,k\g}{1+k\g}{\a}
\end{equation}
for all $n\ge 1,$ as an immediate consequence of \eqref{Mell} - see (3.13) in \cite{GW}. By the  positivity of all parameters in the hypergeometric function, one has $F_{\a,\g}(k) \ge 1$ for all $k\ge 1$ and $\esp[\Y_{\a, \g}^n] \le n!$ for all $n\ge 1,$ so that $\esp[e^{z\Yag}] < \infty$ for all $z < 1$ and the law of $\Yag$ is characterized by its integer moments.

On the other hand, the Euler integral formula yields
$$\frac{F_{\a,\g}(k)}{k}\; =\; \g \int_0^1 x^{k\g -1} (1- \a x)^{-\g}\ dx \; =\; \int_0^1 x^{k -1} (1- \a x^{1/\g})^{-\g}\ dx$$
for all $k\ge 1.$ Setting $\{\Y_{\a,\g,k}, \, k\ge 1\}$ for an independent sequence of random variables with respective densities
$$\frac{k\, x^{k -1} (1- \a x^{1/\g})^{-\g}}{F_{\a,\g}(k)}\,\Un_{(0,1)} (x)$$
and respective expectations
$$\esp[\Y_{\a,\g,k}]\; =\; \frac{k\,F_{\a,\g}(k+1)}{(k+1)\, F_{\a,\g} (k)},$$
the positive martingale
$$\prod_{k=1}^n \lpa \frac{(k+1)\, F_{\a,\g} (k)}{k\,F_{\a,\g}(k+1)}\rpa  \Y_{\a,\g,k}$$
converges a.s. to the positive random variable
$$\prod_{k=1}^\infty \lpa \frac{(k+1)\, F_{\a,\g} (k)}{k\,F_{\a,\g}(k+1)}\rpa  \Y_{\a,\g,k},$$
whose $n$-th integer moment is computed by Fubini's theorem as 
\begin{eqnarray*}
\prod_{k=1}^\infty \lpa \frac{(k+1)\, F_{\a,\g} (k)}{k\,F_{\a,\g}(k+1)}\rpa^{\! n}  \esp[\Y_{\a,\g,k}^n] & = &\prod_{k=1}^\infty \lpa \frac{(k+1)\, F_{\a,\g} (k)}{k\,F_{\a,\g}(k+1)}\rpa^{\! n} \!\lpa\frac{k\,F_{\a,\g} (n+k)}{(n+k)\, F_{\a,\g} (k)}\rpa\\
& = & \lim_{N\to\infty} \prod_{k=1}^N \lpa \frac{(k+1)\, F_{\a,\g} (k)}{k\,F_{\a,\g}(k+1)}\rpa^{\! n} \!\lpa\frac{k\,F_{\a,\g} (n+k)}{(n+k)\, F_{\a,\g} (k)}\rpa\\ 
& = &  \prod_{k=1}^n \, \frac{k\,F_{\a,\g}(1)}{F_{\a,\g}(k)}\,\times\, \lim_{N\to\infty} \lpa\frac{F_{\a,\g} (N+1)\times\cdots\times F_{\a,\g}(N+n)}{\lpa F_{\a,\g}(N+1)\rpa^n}\rpa \\
& = & \lpa F_{\a,\g}(1)\rpa^n \prod_{k=1}^n \, \frac{k}{F_{\a,\g}(k)}
\end{eqnarray*}
for every $n\ge 1,$ where in the last equality we have used the fact that
$$F_{\a,\g}(N+k)\; =\; (1-\a)^{-\g}  \pFq{2}{1}{\g,,1}{1+(N+k)\g}{\frac{\a}{\a -1}}\;\to\; (1-\a)^{-\g}$$
as $N\to\infty$ for every $k=1,\ldots, n,$ by Kummer's identity on the hypergeometric function. By integer moment identification, this shows 
\begin{equation}
\label{prod}
\Yag\; \elaw\; \frac{1}{F_{\a,\g}(1)} \prod_{k=1}^\infty \lpa \frac{(k+1)\, F_{\a,\g} (k)}{k\,F_{\a,\g}(k+1)}\rpa  \Y_{\a,\g,k},
\end{equation}
which leads to
\begin{equation}
\label{produ}
\Xag\;\elaw\;F_{\a,\g}(1)^{\frac{1}{\g}} \prod_{k=1}^\infty \lpa \frac{k\,F_{\a,\g}(k+1)}{(k+1)\, F_{\a,\g} (k)}\rpa^{\frac{1}{\g}} \Y_{\a,\g,k}^{-\frac{1}{\g}}.
\end{equation}
Now for every $k\ge 1,$ the density of the random variable $\Y_{\a,\g,k}^{-\frac{1}{\g}}-1$ reads
$$\frac{k\g}{F_{\a,\g}(k)\, (x+1)^{(k-1)\g +1} (x+1-\a)^{\g}}$$
on $(0,\infty)$ and is hence hyperbolically completely monotone, with the notation of Chapter 5 in \cite{Bo}. By Theorem 5.1.2 in \cite{Bo} this implies that the law of $\Y_{\a,\g,k}^{-\frac{1}{\g}} = 1 + (\Y_{\a,\g,k}^{-\frac{1}{\g}}-1)$ is a GGC for every $k\ge 1$ and by \eqref{produ} and the main result of \cite{Bond} the law of $\Xag,$ too.
\qed
 
\begin{Rem}
\label{bell}
{\em The random variable $\Xag$ has finite negative integer moments of any order and a combination of Theorems 4.1.3 and 4.1.4 in \cite{Bo} shows that the density $\hag$ has a smooth extension on the whole $\rl,$ having set $\hag(x) = 0$ for all $x\le 0.$ With the terminology of Chapters 3 and 4 in \cite{Bo}, this means that the Thorin measure of $\Xag$ is infinite. By  Corollary 1.2 in \cite{Kwas}, this implies that the function $\hag$ is bell-shaped on $(0,\infty),$ viz. for every $n\ge 0$ one has 
$$\sharp\{x > 0\;\slash \;\, h_{\a,\g}^{(n)} (x) =0\}\; = \; n.$$
At the visual level, the bell-shape property means that the function $\hag$ is increasing-then-decreasing ($n=1$) and convex-then-concave-then-convex ($n=2$) on $(0,\infty)$, which is illustrated by all the figures plotted in \cite{GW}. On the other hand, it does not seem possible to derive these basic properties neither directly from \eqref{ID} nor from the non-classical special function representation obtained in \cite{GW}.  
}
\end{Rem}

\subsection{Proof of Theorem \ref{asym}} For every $\lbd \ge 0,$ we define
$$F_{\a,\g}(\lbd)\; =\; \pFq{2}{1}{\g,,\lbd\g}{1+\lbd\g}{\a}\; = \;\lbd\g \int_0^1 x^{\lbd\g -1} (1- \a x)^{-\g}\ dx.$$
Integrating by parts and changing the variable, we obtain
$$F_{\a,\g}(\lbd)\; = \; 1\; +\; \a\g \int_0^1 \frac{1- x^{\lbd\g}}{ (1- \a x)^{\g+1}}\, dx\; = \; 1\; +\; \Phi_{\a,\g}(\lbd),$$
where
$$\Phi_{\a,\g}(\lbd)\; =\; \int_0^\infty (1- e^{-\lbd t})\,\frac{\a\,e^{-\frac{t}{\g}}}{(1- \a e^{-\frac{t}{\g}})^{\g+1}}\, dt$$
is the Laplace exponent of some driftless compound Poisson process $\{N^{\a,\g}_t, \; t\ge 0\}.$ By \eqref{Mom} and (1.2) in \cite{CPY}, we deduce 
$$\esp[\Y_{\a, \g}^n]\; =\; \prod_{k=1}^n \, \frac{k}{1 + \Phi_{\a,\g}(k)}\; =\; \esp[\I_{\a, \g}^n]$$
for every $n\ge 1,$ where 
$$\I_{\a,\g}\: = \; \int_0^\L e^{-N^{\a,\g}_t}\, dt$$
with $\L\sim$ Exp(1) independent from $\{N^{\a,\g}_t, \; t\ge 0\}$. By moment identification, we then have
\begin{equation}
\label{Exp}
\Yag\;\elaw\; \I_{\a,\g}.
\end{equation} 
This alternative representation allows us to apply the recent results in \cite{MS} on the asymptotic behaviour at infinity of the density of the exponential functional of a possibly killed subordinator - see also \cite{H} for the unkilled case. To be more specific, the Laplace exponent of the  killed and driftless compound Poisson process 
$${\tilde N}^{\a,\g}_t\; =\; N^{\a,\g}_t \Un_{\{t \le \L\}} + \infty \Un_{\{t > \L\}}, \qquad t\ge 0,$$ 
is
$$F_{\a,\g}(\lbd)\; = \; (1-\a)^{-\g}\; -\; \a\g \int_0^1 \frac{ x^{\lbd\g}}{ (1- \a x)^{\g+1}}\, dx\; \to \; (1-\a)^{-\g}$$
as $\lbd\to \infty$ by monotone convergence. Moreover, we have
$$F_{\a,\g}'(\lbd)\; = \;\int_0^\infty \frac{\a\g^2\,t\, e^{-t (\lbd\g + 1)}}{(1- \a e^{-t})^{\g+1}}\, dt\; =\; O(\lbd^{-2}), \qquad \lbd\to\infty,$$
and this shows that the assumption (2.6) in \cite{MS} is satisfied. A further elementary asymptotic analysis implies
\begin{eqnarray*}
F_{\a,\g}(\lbd) & = & (1-\a)^{-\g}\lpa 1\, -\, \frac{\a}{(1-\a)\lbd}\int_0^\infty e^{-t(1+ \frac{1}{\lbd\g})} \lpa \frac{1-\a}{1- \a e^{-\frac{t}{\lbd\g}}}\rpa^{\g +1} \!\!dt\rpa\\
& = & (1-\a)^{-\g}\lpa 1\, -\, \frac{\a}{(1-\a)\lbd} \, +\, O\lpa \lbd^{-2}\rpa\rpa, \qquad \lbd\to\infty.
\end{eqnarray*}
Skipping the easy details we then obtain, with the notations of \cite{MS},  
$$\frac{\varphi_*(y)}{y}\; =\; (1-\a)^{-\g}\, -\, \frac{\a}{(1-\a) y} \, +\, O\lpa y^{-2}\rpa\qquad\mbox{and}\qquad \varphi_*'(x)\;\to\; (1-\a)^{-\g}\qquad\mbox{as $x, y\to\infty.$}$$ 
Setting $g_{\a,\g}$ for the density of $\Yag,$ putting everything together and applying Theorem 3.1. in \cite{MS} yields finally the existence of a positive finite constant $c$ such that
\begin{equation}
\label{Asyl}
g_{\a,\g} (x)\; \sim\; c\, x^{\frac{\a}{1-\a}} \, e^{-(1-\a)^{-\g} x}, \qquad x\to\infty.
\end{equation}
The constant $c$ can be identified via the constant $C$ given in the statement of Theorem 3.1 in \cite{MS}, whose expression is unfortunately very complicated. Instead, we will proceed as in \cite{JSW} and use the large integer moments asymptotics. On the one hand, the Laplace approximation entails
$$\frac{\esp[\Y_{\a, \g}^n]}{n!} \; \sim \; \frac{c}{n!}\int_0^\infty x^{n+\frac{\a}{1-\a}} e^{-(1-\a)^{-\g} x}\, dx
\; \sim \; c(1-\a)^{\frac{\g}{1-\a}}\, n^{\frac{\a}{1-\a}} (1-\a)^{\g n}, \qquad n\to\infty.$$  
On the other hand, we have by Kummer's identity on the hypergeometric function
\begin{eqnarray*}
\frac{\esp[\Y_{\a, \g}^n]}{n!}\; =\; \prod_{k=1}^n \, \frac{1}{F_{\a,\g}(k)} & = & n^{\frac{\a}{1-\a}} (1-\a)^{\g n} e^{\frac{\a}{\a-1}\lpa\log n - \sum_{k=1}^n \frac{1}{k}\rpa} \prod_{k=1}^n \, \frac{e^{\frac{\a }{(\a-1)k}}}{G_{\a,\g}(k)}\\
& \sim & \lpa e^{\frac{\a\psi(1)}{\a-1}} \prod_{k=1}^\infty \, \frac{e^{\frac{\a}{(\a-1)k}}}{G_{\a,\g}(k)}\rpa n^{\frac{\a}{1-\a}} (1-\a)^{\g n}, \qquad n\to\infty.
\end{eqnarray*}
Comparing those two estimates gives
$$c \; =\; (1-\a)^{\frac{\g}{\a-1}}\, e^{\frac{\a\psi(1)}{\a-1}} \prod_{k=1}^\infty \, \frac{e^{\frac{\a}{(\a-1)k}}}{G_{\a,\g}(k)}$$
and concludes the proof since $\hag(x) = \g x^{-\g-1} g_{\a,\g}(x^{-\g}).$

\qed

\begin{Rem}
\label{SNLP}
{\em By the binomial theorem, the function
$$t\;\mapsto\; \frac{\a\,e^{-\frac{t}{\g}}}{(1- \a e^{-\frac{t}{\g}})^{\g+1}}$$
is completely monotone. This means that the Bernstein function $F_{\a,\g}(\lbd)$ associated to the killed Poisson process $\{{\tilde N}^{\a,\g}_t, \; t\ge 0\}$ is complete, and it is well-known that this yields the following identification between two complete Bernstein functions
$$\frac{\lbd}{F_{\a,\g}(\lbd)}\; =\; \frac{\Psi(\lbd)}{\lbd}$$
where $\Psi$ is the Laplace exponent of a spectrally negative Lévy process, which drifts towards $\infty$ since $\Psi'(0+) = 1/F_{\a,\g}'(0+) > 0.$ By the same argument as in \cite{JSW}, this shows the further identification
$$\Yag\;\elaw\; \J^{-1}_{\a,\g}$$
where $\J_{\a,\g}$ is the perpetuity of this spectrally negative Lévy processes, whose asymptotic behaviour of the density at zero is given by (2.35) in \cite{PS}, and is of course the same as \eqref{Asyl}. Hence, we do not really need the recent estimates of \cite{H, MS}. We chose this point of view in order to stay with the theme, since we will do need such estimates in the next proof where the connection to spectrally negative Lévy processes is less clear.}
\end{Rem}

\subsection{Proof of Theorem \ref{asymb} (c)} We proceed as before and consider the power transformation $\tYag = \tX_{\a,\g}^\g.$ The following recurrence equation for the Mellin transform ${\hat H}_{\a,\g} (s) = \esp [\tX_{\a, \g}^{s-1}]$ is derived from \eqref{IDB} by Fubini's theorem, exactly as \eqref{Mell} was obtained from \eqref{ID}:
\begin{equation}
\label{Mellb}
{\hat H}_{\a,\g} (s) \; =\; \g\, {\hat H}_{\a,\g} (s+\g)\,\int_\a^1 (u-\a)^\g u^{s-2} \, du.
\end{equation}
Observe that since $\a\in (0,1)$, this shows that ${\hat H}_{\a,\g} (s)$ is finite for every $s\in \rl.$ This also implies the following formula for the positive integer moments of $\tYag:$ 
$$ \esp[\tY_{\a,\g}^n] \; =\; {\hat H}_{\a,\g} (1+\g n) \; =\;\prod_{k=1}^n \, \frac{k}{{\hat \Phi}_{\a,\g}(k)}$$
where, integrating by parts, 
$${\hat\Phi}_{\a,\g}(\lbd)\; =\; \int_0^\infty (1- e^{-\lbd t})\;\,\a\,e^{\frac{t}{\g}}\, (1- \a e^{\frac{t}{\g}})_+^{\g-1}\, dt$$
is the Laplace exponent of some unkilled and driftless compound Poisson process, whose Lévy measure has bounded support. The identification between $\tYag$ and the perpetuity of this Poisson process is then made as in the previous proof. We further decompose 
$${\hat \Phi}_{\a,\g}(\lbd) \; = \; (1-\a)^{\g}\lpa 1\, -\, \frac{\a}{(1-\a)\lbd}\int_0^\infty e^{-t(1- \frac{1}{\lbd\g})} \lpa \frac{1- \a e^{\frac{t}{\lbd\g}}}{1-\a}\rpa_+^{\g -1} \!\!dt\rpa,$$
whose asymptotic analysis is performed as above and leads to the behaviour
$${\hat g}_{\a,\g} (x)\; \sim\; c\, x^{\frac{\a}{1-\a}} \, e^{-(1-\a)^{\g} x}$$
at infinity for the density ${\hat g}_{\a,\g}(x)$ of $\tYag,$ for some positive constant $c,$ which cannot seem to be computed in simple terms of the hypergeometric function, as in Theorem \ref{asym}. 

\qed

\subsection{Proof of Theorem \ref{Momus}} We start with the criterion for $\tX^t_{\a,\g}.$ The argument relies on the following terminal value representation, which is obtained from \eqref{Mellb} as in the proof of Theorem \ref{ggc}:
\begin{equation}
\label{Prodi}
\tYag\; \elaw\; \frac{1}{{\hat \Phi}_{\a,\g}(1)} \prod_{k=1}^\infty \lpa \frac{(k+1)\, {\hat \Phi}_{\a,\g} (k)}{k\,{\hat \Phi}_{\a,\g}(k+1)}\rpa  \tY_{\a,\g,k},
\end{equation}
where the independent random variables $\tY_{\a,\g,k}$ have respective densities
$$\frac{k\, x^{k -1} (1- \a x^{-1/\g})_+^{\g}}{{\hat \Phi}_{\a,\g}(k)}\,\Un_{(0,1)} (x).$$
The easily established log-concave character of the function $x\mapsto (1- \a e^{-x/\g})_+^{\g}$ for every $\a\in(0,1)$ and $\g > 0$ implies by a change of variable that the random variable $\log \tY_{\a,\g,k}$ has a log-concave density for all $k\ge 1.$ By \eqref{Prodi} and the Prékopa-Leindler theorem, this shows that the density of $\log\tYag$ and hence also that of $\log\tXag$ is log-concave, in other words that the function $x \mapsto \thag(e^{tx})$ is log-concave on $\rl$ for all $\a\in (0,1), \g > 0$ and $t\neq 0.$ This means that Condition L on p.11 in \cite{Lin} is satisfied by the density of $\tX^t_{\a,\g}$ for all $t\neq 0.$ For $t > 0,$ Krein's condition, Theorem 10 in \cite{Lin} and a change of variable show that 
$$\tX^{t}_{\a,\g} \;\,\mbox{is $M$-det}\quad\Longleftrightarrow\quad-\int_1^\infty \frac{\log\thag(x^{\frac{2}{t}})}{1+x^2}\, dx \, =\, \infty\quad\Longleftrightarrow\quad t \le 2\g$$
as required, where for the second equivalence we have used Theorem \ref{asymb} (c). For $t < 0,$ we will check the converse Carleman condition. Setting ${\hat \Z}_{\a,\g} = \tX^{-\g}_{\a,\g},$ an iteration using \eqref{Mellb} shows 
$$\esp[{\hat \Z}_{\a,\g}^n]\; =\; {\hat H}_{\a,\g} (1-\g n) \; =\;\prod_{k=1}^n \,\frac{{\hat \Psi}_{\a,\g}(k)}{k}$$
for every $n\ge 1,$ where
$${\hat\Psi}_{\a,\g}(\lbd)\; =\;\int_0^\infty\, (e^{-\frac{t}{\lbd\g}}-\a)_+^\g\, e^t\, dt, \qquad \lbd > 0.$$
We then have the lower bound
$${\hat\Psi}_{\a,\g}(k)\; =\; \a^{(1-k)\g} \int_0^{-k\g\log\a}\!\!\!  (e^{\frac{t}{k\g}}-1)^\g\, e^{-t}\, dt \; \ge\; c\,k^{-\g}\, \a^{-k\g}$$
for some positive constant $c,$ where in the inequality we have used $e^u \ge 1+u$ for all $u\ge 0.$ By means of Hölder's inequality, this yields
$$\esp[{\hat \Z}_{\a,\g}^{tn}]\; \ge\;\esp[{\hat \Z}_{\a,\g}^{[tn]}]^{\frac{tn}{[tn]}}\; \ge\;c^n \lpa [nt]!\rpa^{-\frac{tn(1+\g)}{2[tn]}} \a^{-\frac{tn(1+[tn])}{2}}$$
for every $t > 0,$ which by Theorem 7 and the aforementioned condition L in \cite{Lin} implies that ${\hat \Z}_{\a,\g}^{t}$ is $M-$indet for every  $t > 0,$ as required. We finally proceed to the criterion for $\X_{\a,\g}^{-t}.$ Recall that the case $t < 0$ is irrelevant since $\Xag$ has infinite positive integer moments. Recall also from the proof of Theorem \ref{asym} and Remark \ref{SNLP} that $\Yag = \X^{-\g}_{\a,\g}$ has positive integer moments given by
$$\esp[\Y_{\a,\g}^n] \; =\; \prod_{k=1}^n \, {\tilde \Phi}_{\a,\g}(k)\qquad\mbox{with} \qquad {\tilde \Phi}_{\a,\g}(\lbd)\; =\; \frac{\lbd}{F_{\a,\g}(\lbd)},$$
and that ${\tilde \Phi}_{\a,\g}(\lbd)\sim (1-\a)^\g \lbd$ as $\lbd\to\infty,$ is a complete Bernstein function. With the notation of \cite{MB}, this means that $\Yag$ is a Remainder,  which satisfies Assumption 1 p.581 therein - see the discussion thereafter. By Proposition 2.1, Proposition 2.2, Remark 2.3. and Remark 3.2 in \cite{MB}, we deduce
$$\Y^{t}_{\a,\g} \;\,\mbox{is $M$-det}\quad\Longleftrightarrow\quad-\int_1^\infty \frac{\log g_{\a,\g}(x^{\frac{2}{t}})}{1+x^2}\, dx \, =\, \infty\quad\Longleftrightarrow\quad t \le 2,$$
where the second equivalence follows from \eqref{Asyl}. Changing the variable, this completes the proof.

\qed
 
\begin{Rem}
\label{Mouss}
 {\em (a) The log-concavity of $\thag(e^x)$ on $\rl$ implies that $\thag$ is unimodal and, by the real-analyticity of $\thag$ on $(0,\infty)$ which can be obtained from \eqref{IDB} as in Proposition 3.6 of \cite{GW}, that it is strictly unimodal. The existence of negative moments of arbitrary order, which was observed during the proof of Theorem \ref{asymb} (c), clearly shows that the mode is positive. Observe that $\thag$ is however not necessarily bell-shaped since this property would imply infinite divisibility by Corollary 1.4 in \cite{KS}, and this property is not true for $\g > 1$ as discussed in the introduction after the statement of Theorem \ref{asymb}. We conjecture that $\thag$ is bell-shaped if and only if $\g \le 1.$ As in Remark \ref{bell}, the if part would be a consequence of the GGC property and Theorem \ref{asymb} (b).

\medskip 

(b) The positive integer moments of ${\hat \Z}_{\a,\g}$ can also be expressed as 
$$\esp[{\hat \Z}_{\a,\g}^n]\; =\; m\,\prod_{k=1}^{n-1} \,\frac{{\tilde \Psi}_{\a,\g}(k)}{k}$$
with 
$$m\, =\, \g\int_0^\infty (1-\a e^t)_+^\g \, dt\qquad\mbox{and}\qquad {\tilde \Psi}_{\a,\g}(k)\, =\, \int_0^\infty (e^{\lbd t} -1)\;\,\a\,e^{\frac{t}{\g}}\, (1- \a e^{\frac{t}{\g}})_+^{\g-1}\, dt.$$
Observe that $m = {\tilde \Psi}_{\a,\g}'(0+) >0$ and that ${\tilde \Psi}_{\a,\g}$ is the Laplace exponent of a spectrally positive Lévy process which is precisely the subordinator appearing in the proof of Theorem \ref{asymb} (c). This is compatible with ${\hat \Z}_{\a,\g} = {\hat \Y}_{\a,\g}^{-1}$ and Proposition 2 in \cite{BY} - see also the remark on moment-indeterminacy after Proposition 1 therein. In view of Theorem \ref{asymb} (b), one might ask if $-\log\thag(x)$ would not behave at zero like a power of the logarithm, as is the case for the perpetuity of the standard Poisson process - see Proposition 3 in \cite{BY}.

\medskip

(c) The infinite product representation \eqref{Prodi} involves random variables with support $[\a,1]$ which is bounded away from zero, and whose negative powers are hence never ID. It is not clear whether the product representation \eqref{Prodi} can help establish the presumed ID property for $\tXag$ when $\g\le 1.$} 
\end{Rem}

\subsection{Proof of Theorem \ref{asymb} (a) and (b)} The strict unimodality was discussed in Remark \ref{Mouss} (a), and we hence only need to prove (b). The positivity of the mode implies that $\thag$ is increasing in a neighbourhood of zero and for every $a >0$ and sufficiently small $x,$ we have
$$x^{-a}\thag(x) \;\le\; 2^{1+[a]} x^{1+[a] -a} \int_x^{2x} y^{-1-[a]}\, \thag(y)\, dy\; \to\; 0\qquad\mbox{as $x\to 0,$}$$
by the finiteness of negative integer moments. For the second estimate, we recall that $\tX_{\a,\g}^{-t}$ is $M-$indet and satisfies Condition L in \cite{Lin} for every  $t > 0,$ by the proof of Theorem \ref{Momus}. By a change of variable and Theorem 10 in \cite{Lin}, we have
$$-\int_0^1 \frac{\log\thag(x^t)}{1+x^2}\, dx \; <\; \infty$$
for every $t > 0$, which implies our claim by the same monotonicity argument as above. 
 
\qed

\section{Some further comments}

Comparing \eqref{prod} and \eqref{Exp} gives the identity in law
$$\int_0^\infty e^{-{\tilde N}^{\a,\g}_t}\, dt \;\elaw\; \frac{1}{F_{\a,\g}(1)} \prod_{k=1}^\infty \lpa \frac{(k+1)\, F_{\a,\g} (k)}{k\,F_{\a,\g}(k+1)}\rpa  \Y_{\a,\g,k}$$
where we recall that $\{{\tilde N}^{\a,\g}_t, \; t\ge 0\}$ is the compound Poisson process, killed and driftless, with Laplace exponent
$$\esp\lcr e^{-\lbd {\tilde N}^{\a,\g}_t}\rcr\; =\; e^{-tF_{\a, \g}(\lbd)}, \qquad t, \lbd \ge 0.$$
This identification between the perpetuity of a subordinator and the terminal value of a multiplicative martingale is not a surprise and holds actually in full generality. More precisely, if $\{\sg_t,\; t\ge 0\}$ is a subordinator having Laplace exponent
$$\Phi(\lbd)\; =\; q\; +\; b\lbd\; + \; \int_0^\infty (1-e^{-\lbd t}) \,\pi(dt), \qquad \lbd \ge 0,$$
with $b,q \ge 0$ and $\pi(dt)$ a positive measure integrating $1\wedge t$ on $(0,\infty),$ and if $(q,b,\pi)\neq (0,0,0),$ then the identity
\begin{equation}
\label{Ida}
\int_0^\infty e^{-\sg_t}\, dt \;\elaw\; \frac{1}{\Phi(1)} \prod_{k=1}^\infty \lpa \frac{(k+1)\,\Phi (k)}{k\,\Phi(k+1)}\rpa  \Y_k
\end{equation}
holds, where $\{\Y_k, \;k\ge 1\}$ is an independent sequence of random variables on $(0,1]$ with respective distributions
$$\frac{k}{\Phi(k)} \lpa b\delta_1(dx) \, +\, x^{k-1}(q + \pi(-\log x,\infty))\, \Un_{(0,1)} (x)\, dx\rpa$$
and respective integer moments
$$\esp[\Y_k^n]\; =\; \frac{k\,\Phi(k+n)}{(k+n)\, \Phi(k)},\qquad k,n\ge 1.$$
The identity \eqref{Ida} comes from the fact that the random variables on both sides are characterized by their integer moments given by
$$\prod_{k=1}^n \frac{k}{\Phi(k)}, \qquad n\ge 1,$$
which for the perpetuity follows from (1.2) in \cite{CPY}, and for the terminal value is a consequence of the simple evaluation
$$\frac{\Phi(\lbd)}{\lbd}\; =\; b\; +\; \int_0^1 x^{\lbd-1}\,(q + \pi(-\log x,\infty))\; dx,\qquad \lbd > 0,$$
obtained by an integration by parts and a change of variables, and of the very same reasoning as in the proof of Theorem 1 using
$$\lim_{x\to\infty}\frac{\Phi(x+c)}{\Phi(x)}\; =\; 1$$
for all $c > 0.$ We refer to Theorem 2.22 in \cite{PS} for more general factorizations on perpetuities of Lévy processes, with  different normalization constants. The above details for \eqref{Ida} are given for completeness, and because they are very easy. 

\medskip

Theorems 1 and 2 in the present paper consider the situation $\Y_k = \Y_{\a,\g,k}$ with $q=1, b = 0$ and $\pi(-\log x,\infty) = (1-\a x^{1/\g})^{-\g} -1.$  There are of course many other explicit examples, like the jumpless case $b,q > 0$ and $\pi\equiv 0,$ which gives the factorization
\begin{equation}
\label{beta}
\B_{1, \g}\;\elaw\;\frac{1}{1+\g}\; \prod_{k=1}^\infty \frac{(k+1)(k+\g)}{k(k+1+\g)}\;\Y_k
\end{equation}
where $\g = qb^{-1},$ the $\Y_k$ have respective distributions
$$\frac{k}{k+\g} \lpa \delta_1(dx) \, +\, \g x^{k-1}\, \Un_{(0,1)} (x)\, dx\rpa$$
and, here and throughout, $\B_{a,b}$ stands for the beta distribution with parameters $a,b >0.$ Observe also that size-biasing \eqref{beta} leads to the more general
$$\B_{a, b}\;\elaw\;\frac{a}{a+b}\; \prod_{k=0}^\infty \frac{(k+a+1)(k+a+b)}{(k+a)(k+a+b+1)}\;\Y_k$$
for all $a, b >0,$ where the $\Y_k$ have respective distributions
$$\frac{k + a}{k+a + b} \lpa \delta_1(dx) \, +\, b x^{k+a-1}\, \Un_{(0,1)} (x)\, dx\rpa.$$

\medskip

We conclude this note with a brief focus on another example taken from the recent paper \cite{Ber}, and connected to the local time process $\{{\hat L}_t, \; t\ge 0\}$ at level zero of a noise-reinforced Bessel process with dimension $d\in (0,2).$ Corollary 4.3 in \cite{Ber} relates the random variable ${\hat L}_1$ to the perpetuity ${\hat I}$ of a subordinator given in (4.4) therein. See also \cite{KP} and the references therein for some instances of the same perpetuity as limit laws of one-sided tree destructions. Setting $\a = 1-d/2\in (0,1)$ and $\g = 1- 2p \in (0,\infty)$ where $p$ is the reinforcement parameter, the latter subordinator has characteristics $b = q = 0$ and
$$\pi(-\log x, \infty)\; =\; \frac{(\g/2)^\a \, x^\g (1 - x^{-\g/\a})^{-\a}}{\Ga(\a+1)},$$
with our above notation. Applying \eqref{Ida} leads after some computations to
$${\hat I}\; \elaw\; \esp[{\hat I}]\, \prod_{k=1}^\infty \lpa\frac{\Y_k}{\esp[\Y_k]}\rpa$$
with $\Y_k = \B_{\a(1+\frac{k}{\g}),1-\a}^{\frac{\a}{\g}}$ for all $k\ge 0.$ Setting $\Z_k = \B_{\a(1+\frac{k}{\g}),1-\a}$ for all $k\ge 0,$ Corollary 4.3 in \cite{Ber} reads 
\begin{equation}
\label{LTA}
{\hat L}_1\; \elaw\; \esp[{\hat L}_1]\, \prod_{k=0}^\infty \lpa\frac{\Y_k}{\esp[\Y_k]}\rpa \; \elaw \; \lpa\esp[{\hat L}_1^{\frac{\g}{\a}}]\,\prod_{k=0}^\infty \lpa\frac{\Z_k}{\esp[\Z_k]}\rpa\rpa^{\!\frac{\a}{\g}},
\end{equation}
where the first identity comes from size-biasing as in the above beta example, and the second one from integer moment identification using Fubini's theorem as in the previous section. By the second identity in the Theorem of \cite{JSW} and an easy size-bias analysis whose details are skipped, we obtain finally the simple identity
\begin{equation}
\label{LT}
{\hat L}_1\; \elaw\; \esp[{\hat L}_1]\,\prod_{k=0}^\infty \lpa\frac{k+1}{k+\a}\rpa \B_{\g(1+\frac{k}{\a}),\g(\frac{1}{\a}-1)}.
\end{equation}
In particular, the law of ${\hat L}_1^{-1}$ is a GGC as an independent product of reciprocal beta random variables, which are easily seen to have translated HCM densities as in the proof of Theorem \ref{ggc}.

A relevant random variable is ${\hat \lbd}_1 \elaw {\hat L}_1^{-1/\a},$ since it is the value at one of the inverse local time process 
$${\hat \lbd}_t\; =\;\inf\{s\ge 0, \; {\hat L}_s > t\}, \qquad t\ge 0,$$ 
which is an increasing, $\a-$self-similar, time-homogeneous Feller process on $\rl^+$ - see Corollary 4.2 and (4.2) in \cite{Ber}. It has been shown in Lemma 4.1 of \cite{Ber} that ${\hat \lbd}_t$ is distributed as the power transform of an $\a-$stable subordinator taken at the inverse of an additive functional, but it is not clear to the author whether this gives directly some infinitely divisible properties for ${\hat \lbd}_1$, except in the case without reinforcement $\g = 1$. We have the following partial result.

\begin{Prop}
\label{ggcbis}
The law of ${\hat \lbd}_1$ is a {\em GGC} if $(1-\a)\g \ge \a$ or if $\g\in [1,2].$
\end{Prop}

\proof

For $(1-\a)\g \ge \a,$ we have $b = \g(1/\a-1)\ge 1$ and $1/\a\ge 1$ and the law of 
$$\B_{\g(1+\frac{k}{\a}),\g(\frac{1}{\a}-1)}^{-1/\a}$$ is a GGC for all $k\ge 0$ by Theorem 2 (1) in \cite{BS}, whence the result by \eqref{LT} and the main Theorem of \cite{Bond}. For $\g\in [1,2]$ we use the second factorization in \eqref{LTA} and need to show that the law of 
$$\B_{\a(1+\frac{k}{\g}),1-\a}^{-1/\g}$$ 
is a GGC for all $k\ge 0,$ which is here a consequence of Theorem 2 (3) in \cite{BS}.

\endproof 

It has been conjectured in Section 7 of \cite{Bond} that power transformations of order greater than one leave the GGC property invariant, which would show by the preceding discussion that the law of ${\hat \lbd}_1$ is a GGC for all $\a\in (0,1)$ and $\g > 0.$  However, to the best of our knowledge this problem is still open. Consider finally the renormalized local time
$${\tilde L}_1(d,p)\; =\; \frac{{\hat L}_1}{\esp[{\hat L}_1]},$$
which we reparametrize by $d = 2(1-\a)\in(0,2)$ and $p= (1-\g)/2\in(-\infty,1/2).$ It is easy to see from Theorem 1.2 in \cite{Ber} combined with the method of moments and Stirling's formula that for each fixed $d\in (0,2)$ one has
$${\tilde L}_1(d,p)\;\claw\; 1\quad\mbox{as $p\to-\infty$}\quad\quad\mbox{and}\quad\quad {\tilde L}_1(d,p)\;\claw\; c_\a^{-1}\, \cB(c_\a)\quad\mbox{as $p\to 1/2$}$$
where $c_\a = \a^{\a}/\Ga(1+\a) < 1$ by Gautschi's inequality on the Gamma function, and $\cB(q)$ stands for a Bernoulli random variable with parameter $q\in (0,1).$ The same argument yields 
$${\tilde L}_1(d,p)\;\claw\; 1\quad\mbox{as $d\to 0$}\quad\quad\mbox{and}\quad\quad {\tilde L}_1(d,p)\;\claw\; (1-2p)^{-1}\, \G_{1-2p}\quad\mbox{as $d\to 2$}$$
for each fixed $p\in (-\infty,1/2).$ This shows that the limits of ${\tilde L}_1(d,p)$ as $p\to 1/2$ resp. $d\to 2$ are more dispersed than the limits as $p\to -\infty$ resp. $d\to 0.$ Following \cite{HPRY}, we say that a collection $\{X_t, \, t\in I\}$ of integrable random variables indexed by a real interval $I$ is a peacock, if it is increasing for the convex order viz. $\esp[\psi (X_s)]\le \esp[\psi(X_t)]$ for every $s\le t\in I$ and every convex function $\psi$ such that the expectations exist. The following property is another simple consequence of \eqref{LTA} and \eqref{LT}.

\begin{Prop}
\label{peacock}
One has 

\medskip

{\em (a)} For each $d\in (0,2),$ the family $\{{\tilde L}_1(d,p), \; p < 1/2\}$ is a peacock.

\medskip

{\em (b)} For each $p\in (-\infty,1/2),$ the family $\{{\tilde L}_1(d,p), \; d\in (0,2)\}$ is a peacock.
\end{Prop}

\proof We begin with (a). By the product representation \eqref{LT} with $\g = 1-2p$ and the stability of the convex order by mixtures - see Corollary 3.A.22 in \cite{SSh}, it is enough to show that the mapping  $t\mapsto\B_{ta,tb}$ decreases for the convex order on $(0,\infty)$ for all $a,b > 0.$ Fixing $a,b > 0$ and $s < t\in (0,\infty),$ and setting $f_t$ and $f_s$ for the respective densities of 
$\B_{ta,tb}$ and $\B_{sa,sb}$ on $(0,1)$, an immediate analysis shows that
$$\sharp\{ x\in (0,1),\; f_t(x) = f_s(x)\}\; =\; 2$$
and that $f_s(x) > f_t(x)$ in the neighbourhoods of zero and one. This concludes the proof by Theorem 3.A.44 in \cite{SSh}. The proof of (b) goes along the same lines, but it is more involved since we need to care about supports.  By the first product representation in \eqref{LTA} with $d = 1- \a/2$ and the stability of the convex order by mixtures, we need to show that the mapping
$$\a\; \mapsto\; \T_\a\; =\; \lpa \frac{\Ga(1+ (b+1)\a) \,\Ga(1+(a+b)\a)}{\Ga(1+b\a)\,\Ga(1+ (a+b+1)\a)}\rpa \B_{\a(1+b), 1-\a}^{\a a}$$
decreases for the convex order on $(0,1)$ for all $a,b > 0.$ To do so, we first show that the prefactor giving the right end of the support is a decreasing function in $\a$. Taking the logarithmic derivative and reparametrizing, this amounts to $c\psi(1+c) + f \psi (1+ f) >  d \psi (1+d) + e\psi(1+e)$ for all $0 < c < d,e < f$ with $c+f = d+e,$ where
$$\psi(1+z)\; = \; \psi(1) \, +\, \sum_{n\ge 1} \frac{z}{n(n+z)}$$
is the digamma function, which is a consequence of the strict convexity of $x\mapsto x^2/(n+x)$ on $(0,\infty)$ for all $n\ge 1.$ Setting $t_\a(x)$ for the density of $\T_\a$ over the interior of its support $(0,m_\a),$ we deduce that $t_\a(x) > t_\b (x) = 0$ if $x\in ]m_\b,m_\a[,$ whereas $t_\a(x) = t_\b (x) = 0$ if $x > m_\a$ and $t_\b(x)\to \infty > t_\a(m_\b)$ as $x\uparrow m_\b,$ for all $0 < \a < \b < 1.$ Moreover the densities $t_\a(x)$ and $t_\b(x)$ cross at least once on $]0,m_\b[$ since otherwise we would have $\T_\b \le_{st} \T_\a$ by Theorem 1.A.12 in \cite{SSh}, which is impossible by the equality of expectations. Appealing again to Theorem 3.A.44 in \cite{SSh}, we are reduced to prove that these densities cross exactly once on $]0,m_\b[.$ Evaluating the densities and reparametrizing, this amounts to show that the function $x\,\mapsto\, c_1(1- c_2 x^\delta)^{1/\delta} + x - 1$ crosses the positive axis only once on $(0,1)$ for all $c_1, c_2 < 1 < \delta.$ Computing its second derivative $$-c_1c_2(\delta -1) x^{\delta-2}(1-c_2x^\delta)^{1/\delta -2}\, <\, 0$$ 
shows that the latter function is concave, negative at zero and positive at one, which finishes the proof.

\endproof

\begin{Rem} {\em In the non-reinforced case $p = 0,$ the discussion after Theorem 1.2 in \cite{Ber} recalls the classical fact on local time of recurrent Bessel processes that ${\tilde L_1}$ is a renormalized Mittag-Leffler random variable of index $\a,$ with moment generating function $\esp[e^{z{\tilde L}_1}] = E_\a(\Ga(1+\a) z),$ where 
$$E_\a (z) \; =\; \sum_{n\ge 0} \frac{z^n}{\Ga(1+\a n)}, \qquad z\in\rl,$$
is the classical Mittag-Leffler function. The peacock property of Proposition \ref{peacock} (b) implies that for every $z\in\rl,$ the mapping 
$$\a\, \mapsto\, E_\a\lpa(\Ga(1+\a)\, z\rpa$$ 
decreases on $(0,1).$ This had been observed for $\a\in [1/2,1)$ in our previous paper \cite{TSEJP} as a consequence of Theorem B therein. This monotonicity property specifies the way the Mittag-Leffler function interpolates between the hyperbolic curve $(1-z)_+^{-1}$ at the limit $\a \to 0$ and the exponential curve $e^z$ at the limit $\a\to 1.$}

\end{Rem}


\begin{thebibliography}{10}

\bibitem{Ber}
J.~Bertoin. On the local times of noise reinforced Bessel processes. {\em Ann. Henri Lebesgue} {\bf 5}, 1277-1294, 2022.

\bibitem{BY}
J.~Bertoin and M.~Yor. On the entire moments of self-similar Markov processes and exponential functionals of Lévy processes. {\em Ann. Fac. Sci. Toulouse} {\bf 11} (1), 33-45, 2002.

\bibitem{Bo}
L.~Bondesson. {\em Generalized Gamma convolutions and related classes of distributions and densities.} Lect. Notes
Stat. {\bf 76}, Springer-Verlag, New York, 1992.

\bibitem{Bond}
L.~Bondesson. A class of probability distributions that is closed with respect to addition as well as multiplication of independent random variables. {\em J. Theor. Probab.} {\bf 28} (3), 1063-1081, 2015.

\bibitem{BS}
P.~Bosch and T.~Simon. On the infinite divisibility of inverse Beta distributions. {\em Bernoulli} {\bf 21} (4), 2552-2568, 2015.

\bibitem{CPY}
P. Carmona, F. Petit and M. Yor. On the distribution and asymptotic results for exponential functionals of Lévy processes. In:  Exponential functionals and principal values related to Brownian motion. {\em Rev. Mat. Iberoamericana}, 73-121, 1997.

\bibitem{GH} 
P.~E.~Greenwood and G.~Hooghiemstra. On the domain of attraction of an operator between supremum and sum. {\em Probab. Theory Related Fields} {\bf 89} (2), 201-210, 1991. 

\bibitem{H}
B.~Haas. Precise asymptotics for the density and the upper tail of
exponential functionals of subordinators. Available at {\tt arXiv:2106.08691.}    

\bibitem{HPRY}
F.~Hirsch, C.~Profeta, B.~Roynette and M.~Yor. {\em Peacocks and Associated Martingales, with Explicit Constructions.} Bocconi \& Springer Series {\bf 3}, Springer-Verlag, Milan, 2011.

\bibitem{GHB} 
G.~Hooghiemstra and P.~E.~Greenwood. The domain of attraction of the
$\a$-Sun operator for type II and type III distributions. {\em Bernoulli}, {\bf 3} (4), 479-489, 1997.
 
\bibitem{JSW}
W.~Jedidi, T.~Simon and M.~Wang. Density-solutions to a class of integro-differential equations. {\em J. Math. Anal. Appl.} {\bf 458}, 134-152, 2018.

\bibitem{KP}
M.~Kuba and A.~Panholzer. A note on the limit law of one-sided tree destruction. Available at {\tt arXiv:2301.04025.} 

\bibitem{Kwas}
M.~Kwa\'snicki. A new class of bell-shaped functions. {\em Trans. Amer. Math. Soc.} {\bf 373} (4), 2255-2280, 2020.

\bibitem{KS}
M.~Kwa\'snicki and T.~Simon. Characterisation of the class of bell-shaped functions. {\em Math. Z.} {\bf 301}, 2659-2683, 2022.

\bibitem{Lin}
G.-D.~Lin. Recent developments on the moment problem. {\em J. Statist. Dist. Appl.} {\bf 4}, Paper 5 (17 pages), 2017.

\bibitem{MS} 
M.~Minchev and M.~Savov. Asymptotic of densities of exponential functionals of subordinators. Available at {\tt arXiv:2104.05381.}  

\bibitem{PS} 
P.~Patie and M.~Savov. Bernstein-gamma functions and exponential functionals of Lévy processes. {\em Electron. J. Probab.} {\bf 23} (75), 1-101, 2018.

\bibitem{SSh}
M.~Shaked and J.~G.~Shanthikumar. {\em Stochastic orders and their applications.} Springer Verlag, New York, 2007.

\bibitem{MB}
T.~Simon. Moment problems related to Bernstein functions. {\em Ann. Fac. Sci. Toulouse} {\bf 29} (3), 577-594, 2020.

\bibitem{TSEJP}
T.~Simon. Comparing Fr\'echet and positive stable laws. {\em Electron. J. Probab.} {\bf 19} (16), 1-25, 2014. 

\bibitem{GW} 
N.~S.~Witte and P.~E.~Greenwood. On the Density arising from the Domain of Attraction of an Operator between Sum and Supremum: the $\alpha$-Sun Operator. Available at {\tt arXiv:2011.14455.} 

\end{thebibliography}
\end{document}